\documentclass[12pt]{article}

\begin{document}

\begin{center}

About equality of the cosine transform to the sine transform of
Fourier and  the transform of Laplace


%
 \end{center}
\begin{center}
Andrei Pavlov-Maxorin V.

 Moscow, Russia,MIREA-RTU, Vernadskogo
avenu , 78, higher mathematics-1.

login11@umail.ru
 \end{center}

 44A10, 42A38



\begin{center}
Abstract
\end{center}

 Regularity
 of the  transform of Laplace in the opened area of 0 is
 proved with help of some methods of the transform of Fourier. The class of the transform of Laplace
 from the transform of Fourier
  is considered from some functions  without a regularity in null.
  The functions are regular in the opened area of 0.

%





%


\section{Introduction}

An article is by continuation of the \cite{11,22} articles ,
partly published in the \cite{3,4,5} works.

The main part of the  \cite{11,22,3,4,5} works is considered in
the proposition 1. As a result we get the  remarks 1,2   too.


In the next part of article  we consider the regularity of the
double transform of Laplace ( the  theorem 1, the remarks 1, 2).
The theorem proved in this part have a general-mathematical
character and are easily checked up with help of the proposition 1
.

 With help of the   theorem 1 and the remark 2   it is simply to prove a
some  theorems related  with  the transform of Fourier and Laplace
\cite{11,22,3,4,5} ( for instance, about the  inverse operator of
the transform of Laplace, using only positive values of  the
transform of Laplace on the $[0,+\infty)$ \cite{5}). The  theorems
are not by the theme of the the article and require the separate
study in connection with the remarks 1,2. Some  results in the
direction were formulated in the works  \cite{3,4,5}.

The fact about double decomposition on the elementary
 fractions is considered in last part of the article. In opinion of author the
  fact  underlines some interest to the remark 2.

By definition,
$$ L_{\pm} Z(t)(x)=\int\limits_0^{\infty}e^{\pm
xt}Z(t)dt,x\in[0,\infty),$$ (we will use $L_{\pm}
Z(t)(\cdot)(x)=L_{\pm} Z(t)(x)$ too),
$$ F_{\pm} u(t)(\cdot)(p)=\int\limits_{-\infty}^{\infty}e^{\pm
pit}u(t)dt,p\in(-\infty,\infty), L_{+}=L
 ,$$
$$ C^0 u(t)(\cdot)(x)
   =\int\limits_0^{\infty}\cos{xt}u(t)dt, S^0 u(t)(\cdot)(x)
   =\int\limits_0^{\infty}\sin{xt}u(t)dt,x\in(-\infty,\infty),$$

$$
F^{0}_{\pm} u(t)(\cdot)(p)=\int\limits_{0}^{\infty}e^{\pm
pit}u(t)dt,p\in(-\infty,\infty).$$

\section {The regularity of Laplace  transform in $|z|<a>0$}.

In the section we use the Y1 condition.

The {\bf Y1 condition} takes place for the $u(p)$ function, if the
$u(p)$ function is regular  for all p without only k points
$z_1,\ldots,z_k,z_j\notin(-\infty,\infty)\bigcup(-i\infty,i\infty),k=0,1,\ldots$,
$u(0)=0$,
  and
$$\max[|u(p)|,|du(p)/dp|,|d^2u
 (p)/p^2|]|p^{2+\delta}|\rightarrow0,|p|\to\infty,$$
 $\delta>0$,
 $\delta=const.
 $

 We use the Ch1 condition too.

{\bf Ch1 condition.}

 The $u(p)$ function is regular in $K_{++}=\{p:Im\,p\ge 0\}\bigcap
\{p:Re\,p\ge0\}$ or in   $K_{+-}=\{p:Im\,p\le 0\}\bigcap
\{p:Re\,p\ge 0\}$.

{\bf Proposition 1.}

The $l_1(p)+l_2(p)=L_1(p)+L_2(p)=Q(p)$ equality takes place  for
$p=it\in\{p:|p|< \varepsilon>0\}$ as for $t>0$ so as for $t<0$,
$t\in(-\infty,\infty)$, if
$$l_1(p)+l_2(p)=L_1(p)+L_2(p)=Q(p),p=iy\in\{p:|p|< \varepsilon>0\},y\in(0+\infty),$$
 if the $Q(p)$ function
 is regular in $|p|<\varepsilon>0$, if
the  $l_1(p),L_1(p)$ functions are regular in left part of the
complex plane,and the $l_2(p),L_2(p)$ functions are regular in
right part of the plane, if the functions are continuous on the
$(-i\infty,i\infty)$ line, ( for the $u(p)$
  function
  the 
  Y1 condition
 takes place  $Re\,u(t)=u(t),t\in[0,\infty)$).

%

{\bf Proof.}

The $l_1(p)+l_2(p)$ values for $p=it,t>0$  are equal to the
 $L_1(p)+L_2(p)$ values in the right part of the plane, and we
 obtain, that the values are equal as to the
 $L_1(it)+L_2(it),t<0,$ values so as to the
$l_1(it)+l_2(it),t<0,$ values  ( as the result of the
$Q(p)=l_1(p)+l_2(p),p=it,$ function in the left and right parts of
the plane on the $p=it$).

The proposition 1 is proved.

{\bf Remarks 1.}

%

 The
$LLu(x) (\cdot)(z)$ function is regular in
  $\{z:|z|<\varepsilon\}$ for some $\varepsilon>0$.
  if for the  $u(p)$ function the Y1 condition takes place.

{\bf Proof.}

For the
$u_{-}(t)=u(t),t\in[0,+\infty),u_{-}(t)=-u(t),t\in(-\infty,0),$
expression we obtain
$$Q(y)=2\pi
u(y)=2Re\,l_1(x)=l_1(p)+l_2(p)=L_1(p)+L_2(p)=2Re\,L_1(y),p=iy,$$
$y\in[0,+\infty),$ where
$$l_1(p)=LF_{+} u(t)(\cdot)(p), l_2(p)=L_{+}F_{-}
u(t)(\cdot)(p),$$
$$L_1(p)=LF_{+} u_{-}(t)(\cdot)(p), L_2(p)=L_{+}F_{-}
u_{-}(t)(\cdot)(p).$$

From the proposition 1 we get $l_1(p)+l_2(p)=
L_1(p)+L_2(p),p=it,t\in(-\infty,\infty)$, and $l_1(p)-L_1(p)=
L_2(p)-l_2(p),p=it,t\in(-\infty,\infty)$ for all p. We obtain
$l_1(p)-L_1(p)= L_2(p)-l_2(p)\equiv0$ (with help of the Y1
condition). As result the $L_1(p)=LF_{+}
u_{-}(t)(\cdot)(p),L_2(p)$ functions are regular in
$|p|<\varepsilon>0$.

 From $l_1(p)+L_1(p)=2LF^0_{+}
u(t)(\cdot)(p)=2(1/i)LLu(t)(\cdot)(ip)$ we get the remark 1.

{\bf Remarks 2.}

$$0\equiv l_1(ix)+L_1(ix)=
\pm |S^0C^0 u(t)(\cdot)(x)+C^0S^0u(t)(\cdot)(x))|\equiv
0,x\in[0,+\infty),$$ if for the $u(p)$
  function
  the 
  Y1 condition
 takes place, and $Re\,u(t)=u(t),t\in[0,\infty)$.
%
%
%


%

%

{\bf Theorem 1}

 The functions
$$LF^0_{+} u(x)(\cdot)(z)= \int\limits_{0}^{\infty}e^{
-zt} dt\int\limits_{0}^{\infty}e^{ itx}u(x)dx=i LLu(x)
(\cdot)(iz),\,\, LLu(x) (\cdot)(z),
$$  are regular in the
 area  $\{z:|z|<\varepsilon\}$ for some $\varepsilon>0$,
  if for the $u(p)$
  function
  the 
  Y1 condition
 takes place, and $Re\,u(t)=u(t),t\in[0,\infty)$.


{\bf Proof}

We can use the proposition 2 \cite{3,4,5}.

{\bf Proposition 2.}

The equalities $$LF^0_{+} u(x)(\cdot)(v)=iF^0_{-} L
u(x)(\cdot)(v),v\in[0,\infty),$$
$$LC^0u(x)(\cdot)(v)=S^0Lu(x)(\cdot)(v),
LS^0u(x)(\cdot)(v)=C^0Lu(x)(\cdot)(v),v\in[0,\infty)$$  take
place, if for the $u(p)$ function the Y1 condition takes place
(The similar equality $LF^0_{-} u(x)(\cdot)(v)=-iF^0_{+} L
u(x)(\cdot)(v),v\in[0,\infty)$ takes place too).

{\bf Proof.}

We get the first formula  after the change of order of integration
in both parts of the first equality,(if $u(0)=0$, it is obviously
with help of the expressions
$$|F^0_{-} u(x)(\cdot)(t)|\le|(du(x)/dx|_{x=0})/t^2|+|(1/t^2)F^0_{-}
(d^2u(x)/dx^2)(\cdot)(t)|\le c_1/t^2,t\to\infty,$$   $
 c_1=const.,$
$c_1<\infty $, \cite{6}).


With help of the proposition 2 we obtain, that the $F^0_{-} L
u(x)(\cdot)(p)=l_{-}(p)$ function is defined for all $Im\,p<0$,
and $$\lim_{p\to iy,Im\,p<0}F^0_{-} L u(x)(\cdot)(p)=F^0_{-} L
u(x)(\cdot)(iy),y\in(-\infty,\infty),$$( it is obviously,   if
$u(0)=0$ ; as in the proposition 2 we use the formula of
integration on parts \cite{6}).

Similar facts take place for a similar function $l_{+}(p)=F^0_{+}
L u(x)(\cdot)(p)$; the  function is definite from other side of
plane $Im\,p>0$.

We suppose $u(-p)=-u(p)$.

We can write
 $$F^0_{-} L u(x)(\cdot)(p)+F^0_{+} L
u(x)(\cdot)(p)=2C^0Lu(x)(\cdot)(p)=F(p),p=y,y\in[0,\infty),$$ if
$u(-p)=-u(p)$, or
$$l_{-}(p)+l_{+}(p)=F(p),$$where $F(p)$ are regular in
$\{p:|Im\,p|<A\} \bigcup \{p:|Re\,p|<A\}$ for some $A>0$, if the
$u(p)$ function is regular as in the Y1 condition (the fact is
well-known \cite{2,5,6}).

To prove the fact  for $p=y\in (-\infty,0]$ we can define a new
functions $l^{to+}_{-}(p),l^{to-}_{+}(p)$:
$$l^{to+}_{-}(p)=l_{-}(p),Im\,p\le 0,$$ $$l^{to-}_{+}=l_{+}(p),
Im\,p\ge 0,$$where $l^{to+}_{-}(p)$ is an analytical continuation
of the $l_{-}(p),Im\,p\le 0$ function from  the lower part of
plane to the overhead part of plane $\{p:Im\,p\ge 0\}$;
$l_{+}^{to-}(p)$ is an analytical continuation of the
$l_{+}(p),Im\,p\ge 0$ function from from  the overhead part of
plane to the lower part of plane $\{p:Im\,p\le 0\}$ \cite{2}.

The equality $l^{to+}_{-}(p)+l_{+}(p)=F(p),Im\,p\ge 0$ repeats the
main  equality $l_{-}(p)+l_{+}(p)=F(p)$, but in the
$\{p:Im\,p\ge0\}$ area; the  equality
$l_{-}(p)+l^{to-}_{+}(p)=F(p)$ repeats the main  equality
$l_{-}(p)+l_{+}(p)=F(p)$, but in the $\{p:Im\,p\le0\}$ area, where
the $F(p)$ function is regular in $\{p:|Im\,p|<A\}
\bigcup\{p:|Re\,p|<A\}$ for the $u(-p)=-u(p)$ function \cite{2}.
We obtain, that
$$l^{to+}_{-}(p)+l_{+}(p)=l_{-}(p)+l^{to-}_{+}(p),p=y\in[0,\infty).$$

But the same equality takes place
 and for the  $p=y\in (-\infty,0]$ (we use, that both functions
$l^{to+}_{-}(p)+l_{+}(p)=F(p),l_{-}(p)+l^{to-}_{+}(p)=F(p)$ are
equal to the regular $F(p)$ function in different parts of the
plane for the $u(-p)=-u(p)$ function \cite{2}).

We get
$$l^{to+}_{-}(p)+l_{+}(p)=l_{-}(p)+l^{to-}_{+}(p),p=y\in(-\infty,\infty).$$

 The   functions
$l^{to+}_{-}(p),l_{+}(p),l_{-}(p),l^{to-}_{+}(p)$ are the
transforms of Laplace in area of definition, and the functions
 are regular in area of definition \cite{2} (from the proposition 1,2) with  values on the boundary. The
  $l^{to+}_{-}(p)
,l^{to-}_{+}(p))$ functions are regular in the area of regularity
of the sums from the the lemma 1.

{\bf Lemma 1.}


1. The  function $LF^0_{+} u(x)(\cdot)(p)=iF^0_{-} L
u(x)(\cdot)(p)=il_{-}(p)=il^{to+}_{-}(p)$ is regular for all
$\{p:Im\,p\ge 0\}$ (we consider the branch of the function,
passing through $\{p:p=iy,y=Im\,p>0\}$ and $\{p:p=x,x=Re\,p>0\}$,
where the $F^0_{-} L u(x)(\cdot)(p))$ values not defined), if for
the $u(p)$ function or for  the  $V_i(p)$ functions  the Y1,Ch1
conditions take place $i=1,2$,  $u(p)=V_1(p)+V_2(p)$  (see the
part 2 of the  remark  4 too).


2. The function  $LF^0_{-} u(x)(\cdot)(p)=-iF^0_{+} L
u(x)(\cdot)(p)=-il_{+}(p)=l^{to-}_{+}(p)$ is regular  for all
$\{p:Im\,p\le 0\}$ (we consider the branch of the function,
passing through $\{p:p=iy,y=Im\,p<0\}$ and $\{p:p=x,x=Re\,p>0\}$,
where the $F^0_{+} L u(x)(\cdot)(p))$ values not defined), if for
the $u(p)$ function  or for  the  $V_i(p)$ functions  the Y1,Ch1
conditions take place $i=1,2$,  $u(p)=V_1(p)+V_2(p)$ (see the part
2 of the  remark  4 too).

{\bf Proof.}
If the $u(p)$ function is regular in $K_{++}=\{p:Im\,p\ge
0\}\bigcap \{p:Re\,p\ge0\}$,  after integration along the  line
$L_{+}=L_{1}\bigcup L_{2}\bigcup L_{3}$ of the  $(1/ip+z)  u (z) $
function anticlockwise , $L_{1}=[0,R]$,
$L_{2}=\{p:p=Re^{i\varphi},0\le\varphi\le\pi/2\}$, $L_{3}=[iR,0]$,
we obtain, that


$$l^{to+}_{-}(p)=l_{-}(p)=F^0_{-} L u(x)(\cdot)(p)=
\int\limits_{0}^{\infty}(1/ip+x)  u (x) dx=$$
$$=
\int\limits_{0}^{\infty}(1/ip+ix_1)  u (ix_1) dix_1= (1/i)L L
u(ix_1)(\cdot)(p),p\in(0,+\infty),$$  as for
$\{p:Re\,p\notin(-\infty,0)\}$ so as  for all $\{p:Im\,p\ge0\}$.
We use the Y1,Ch1 conditions (the $u(p)$ function is regular in
$K_{++}=\{p:Im\,p\ge 0\}\bigcap \{p:Re\,p\ge0\}$; for $Im\,p=0$ we
use the $u(0)=0$ condition  for proof of continuity on the
$(-\infty,\infty) $ axis with help of the proposition 1 \cite{6}).

We obtain, that the such sum is regular for all $\{p:Im\,p\ge
0,Re\,ip\le0\}$, and the function
 $$\int\limits_{0}^{\infty}(1/ip+ix_1)  u (ix_1) dix_1=(1/i)L L
u(x)(\cdot)(p),Re\,p\notin(-\infty,0)$$  is regular in
$\{p:Im\,p\ge 0,Re\,ip\le0\}$ with the values on the boundary line
$\{p:Im\,p= 0,Re\,p\le 0\}$
   (with help  of the formula of integration on parts
  as in proposition 2  \cite{2,6}).

   For the function $F^0_{+} L
u(x)(\cdot)(p)=l_{+}(p)=l^{to-}_{+}(p),Im\,p<0$  we use
$l^{to-}_{+}(p)=\overline{l^{to+}_{-}(\overline{p})}$, as the
branches of the functions
$l_{+}(p)=\overline{l_{-}(\overline{p})}$
 \cite{2}
 (with help of  the formula $F^0_{+}
L u(x)(\cdot)(p)=\overline{F^0_{-} L u(x)(\cdot)(\overline{p})}$
 on the $[0,+\infty)$ line) by the theorem of Riemann about
the analytical continuation across the $(-\infty,\infty)$ line
\cite{2}, and the function $l_{+}(p)=
l^{to-}_{+}(p)=\overline{l^{to+}_{-}(p)(\overline{p})}$ is defined
and regular in $\{p:Im\,p\le0\}$  (or for the $u(p)=V_1(p)+V_2(p)$
function in connection with the part 2 of the  remark  4 too) .

If the $u(p)$ function is regular in $K_{+-}=\{p:Im\,p\le
0\}\bigcap \{p:Re\,p\ge0\}$,  after integration along the line
$L_{+}=L_{1}\bigcup L_{2*}\bigcup L_{3*}$ of the  $(1/ip+z)  u (z)
$ function anticlockwise, $L_{1}=[0,R]$,
$L_{2*}=\{p:p=Re^{i\varphi},-\pi/2\le\varphi\le0\}$,
$L_{3*}=[-iR,0]$, we obtain, that
$$F^0_{+} L u(x)(\cdot)(p)=l_{+}(p)=
-\int\limits_{0}^{\infty}(1/ip-x)  u (x) dx=$$
$$-\int\limits_{0}^{-\infty}(1/ip-ix_1)  u (ix_1) dix_1 =  $$

$$=
\int\limits_{0}^{+\infty}(1/ip+ix_2)  u (-ix_2) dix_2= L L
u(-ix_1)(\cdot)(p),p\in(0,+\infty).$$

The further proof of lemma repeats  proof of the first part (we
use, that $Re\,u(t)=u(t),t\in[0,\infty)$, but
$u(p)=V_1(p)+V_2(p)$).

 With help  of  lemma 1 we get
$$|l^{to+}_{-}(p)+l_{+}(p)|\le C=const.,Im\,p\ge
0;\,\,l_{-}(p)+l^{to-}_{+}(p)\le C=const.,Im\,p\le 0,$$
$C<\infty$.  Both sums $l^{to+}_{-}(p)+l_{+}(p)$,
$l_{-}(p)+l^{to-}_{+}(p)$ are regular in in area of definition and
continuous on the boundary \cite{2,6}.

We proved $$l^{to+}_{-}(p)+l_{+}(p)= C,C=const.,C<\infty$$
 for all
p  \cite{2}, and
$2C^0Lu(x)(\cdot)(p)=F(p)=l^{to+}_{-}(p)+l_{+}(p)\equiv
C_1=const.,C_1<\infty$ \cite{2,6} or $l^{to+}_{-}(p)=-l_{+}(p)$
for all p including $p\in(-\infty,\infty)$.

We obtain, that the $F(p)$ function is regular \cite{2} with the
values $2C^0Lu(x)(\cdot)(p)=F(p),p\in(-\infty,\infty)$.

 We can use, that the function $LF^0_{+}
u(x)(\cdot)(p)=iF^0_{-} L u(x)(\cdot)(p)$  is regular in
$\{p:|Re\,p|<\varepsilon\}\bigcup \{p:|Im\,p|<\varepsilon\}$ for
some $\varepsilon>0$(it is well-known fact \cite{2,5,6} in Y1
condition for the u(p) function).


We proved, that
$$F^0_{-} L u(x)(\cdot)(p)=F(p)-F^0_{+} L
u(x)(\cdot)(p)$$ is the analytical continuation from
 from one side of plane on other \cite{2}.

 For the $u(-p)=u(p)$ we can use
$$F^0_{-} L u(x)(\cdot)(p)-F^0_{+} L
u(x)(\cdot)(p)=2iS^0Lu(x)(\cdot)(p)=F(p),p=y,y\in[0,\infty),$$
further by analogy with the first part ( with help of the lemma 1
and the $u(0)=0$ condition).

 The theorem 1 is proved.

 From the theorem 1 we obtain the  remark 3.

{\bf Remark  3.}

 The
theorem 1 takes place for the functions
$u(p)=v(p)+v(-p),u(p)=v(p)-v(-p)$, if all the essential points
\cite{2} of the $v(p)$ function are placed or in $K_{++}$ ( or all
the essential points are  in $K_{-+}$).

From the remark 1 we get the first part of the remark 2.


The second part is easily proved by the methods of  the work
\cite{1,6}.

 {\bf Remark  4.}

 1.

The theorem 1 takes place for the function
$u(p)=(v_1(p)-v_1(-p))+(v_2(p)-v_2(-p))$, if the essential points
\cite{2} of the $v_i(p)$ function are placed or in $K_{++}$  or
 in $K_{-+}$, $i=1,2$, and for the
$u(p)$ functions  the conditions of the theorem 1 take place, .

2.

 The   function $u(-p)=-u(p)$ can be presented in the form $u(p)=(v_1(p)-v_1(-p))+(v_2(p)-v_2(-p))$,
 if for the $u(p)$ function the Y1 condition takes place.

We will get a similar result, if to apply the $ z_1=1/z$ inversion
and the $w=e^{i\varphi}z_1,\varphi=\pi/4,$ function; we use, that
 the same result we obtain for  the  functions in the reverse order with $\varphi=-\pi/2$ (by analogy
 for $-\varphi$).


\section{Conclusion}


We will mark the fact  about double decomposition on the
elementary
 fractions:
 $$p[1/(p-1)-1/(p+1)]=1/(p-1)+1/(p+1),
 p/(p-1)^2-p/(p^2-1)=1/(p-1)^2+1/(p^2-1),$$ $p\ne 1,-1.$
 The  fact  in opinion of author underlines interest to the theorem 1 and  to the consequences
  of the  theorem 1.





\end{document}